\documentclass[letterpaper, 10 pt, conference]{ieeeconf}  

\IEEEoverridecommandlockouts                              
\usepackage{graphicx}

\usepackage{color}
\usepackage{array}
\usepackage{phdalgo}
\usepackage{xspace} 
\usepackage{hyperref} 

\usepackage{mathtools}

\usepackage{verbatim}
\usepackage{cite}
\usepackage{paralist}
\usepackage{tikz}
\usetikzlibrary{shapes,arrows,automata,decorations,trees, backgrounds, patterns, calc}

\newcommand{\N}{\mathbb{N}}					
\newcommand{\R}{\mathbb{R}}					

\newcommand{\cB}{\mathcal{B}}
\newcommand{\cF}{\mathcal{F}}
\newcommand{\etc}{\textit{etc}}
\newcommand{\ie}{\textit{i.e.}}
\newcommand{\xopt}[7]{$x_{#1} := (#2, #3, #4, #5, #6, #7)$}
\newcommand{\fsa}{f_{\text{sa}}}
\newcommand{\fpop}{f_{\text{pop}}}
\newcommand{\Kpop}{K_{\text{pop}}}
\newcommand{\Xopt}{x_\mathit{opt}}
\newcommand{\Xsdp}{x_\mathit{sdp}}

\newcommand{\Ksa}{K_{\mbox{\scriptsize sa}}}
\newcommand{\Df}{\mathcal{D}(f)}
\newcommand{\Hf}{\mathcal{D}^2(f)}
\usepackage{amsmath}
\usepackage{amsfonts}
\usepackage{amssymb} 
\usepackage{latexsym}
\usepackage{multirow}
\numberwithin{equation}{section}

\DeclareMathOperator{\parab}{par}

\renewcommand{\geq}{\geqslant}
\renewcommand{\leq}{\leqslant}

\newtheorem{theorem}{Theorem}[section]
\newtheorem{lemma}[theorem]{Lemma}

\newtheorem{definition}[theorem]{Definition}

\newtheorem{example}[theorem]{Example}

\usepackage{xcolor}
\definecolor{darkgreen}{rgb}{0.0 0.5 0.0}
\definecolor{whitegreen}{rgb}{0.0 0.75 0.0}
\definecolor{whiteblue}{rgb}{0.0 0.0 1.2}
\definecolor{darkred}{rgb}{0.8 0.0 0.0}

\title{\LARGE \bf
Certification of inequalities involving transcendental functions: combining SDP and max-plus approximation}

\author{Xavier ALLAMIGEON$^{1}$, St{\'e}phane GAUBERT$^{2}$, Victor MAGRON$^{3}$ and Benjamin WERNER$^{4}$
\thanks{The research leading to these results has received funding from the European Union's $7^{\text{th}}$ Framework Programme under grant agreement nr. 243847 (ForMath).}
\thanks{$^{1}$INRIA and CMAP, \'Ecole Polytechnique, Palaiseau, France,  
        {\tt\footnotesize Xavier.Allamigeon at inria.fr}}%
\thanks{$^{2}$INRIA and CMAP, \'Ecole Polytechnique, Palaiseau, France,         
        {\tt\footnotesize Stephane.Gaubert at inria.fr}}%
\thanks{$^{3}$INRIA and LIX, \'Ecole Polytechnique, Palaiseau, France,      
        {\tt\footnotesize magron at lix.polytechnique.fr}}%
\thanks{$^{4}$INRIA  and LIX, \'Ecole Polytechnique, Palaiseau, France, 
        {\tt\footnotesize benjamin.werner at polytechnique.edu}}%
\thanks{The first author was partially supported by the ForMath EU FP7 STREP FET project, number 243847}
\thanks{Appears in the Proceedings of the European Control Conference ECC'13, July 17-19, 2013, Zurich, pp. 2244--2250, \copyright\ EUCA 2013}
}

\begin{document}

\maketitle
\thispagestyle{empty}
\pagestyle{empty}

\begin{abstract}

  We consider the problem of certifying an inequality of the form
  $f(x)\geq 0$, $\forall x\in K$, where $f$ is a multivariate
  transcendental function, and $K$ is a compact semialgebraic set. We
  introduce a certification method, combining semialgebraic
  optimization and max-plus approximation. We assume that $f$ is given
  by a syntaxic tree, the constituents of which involve semialgebraic
  operations as well as some transcendental functions like $\cos$,
  $\sin$, $\exp$, etc. We bound some of these constituents by suprema
  or infima of quadratic forms (max-plus approximation method,
  initially introduced in optimal control), leading to semialgebraic
  optimization problems which we solve by semidefinite
  relaxations. The max-plus approximation is iteratively refined and
  combined with branch and bound techniques to reduce the relaxation
  gap.  Illustrative examples of application of this algorithm are
  provided, explaining how we solved tight inequalities issued from
  the Flyspeck project (one of the main purposes of which is to
  certify numerical inequalities used in the proof of the Kepler
  conjecture by Thomas Hales).

\end{abstract}

\section{INTRODUCTION}
\label{sect:intro}
\subsubsection*{Inequalities involving transcendental and semialgebraic functions}
Given a multivariate transcendental real function $f : \R^n \to \R$
and a compact semialgebraic set $K \in \R^n$, we consider the
following optimization problem:
\begin{equation}
\label{eq:f}
f^*  :=  \inf_{x \in K} f (x) \enspace,
\end{equation}

The goal is to find the global minimum $f^*$ and a global minimizer
$x^*$. We shall also search for certificates to assess that:
\begin{equation}
\label{eq:ineq}
\forall x \in K, \quad f(x) \geq 0 \enspace.
\end{equation}

A special case of Problem~\eqref{eq:f} is semialgebraic
optimization. Then, $f = \fsa$ belongs to the algebra $\mathcal{A}$ of
semialgebraic functions
which extends multivariate polynomials\linebreak by allowing arbitrary
composition of
$(\cdot)^{p}$, $(\cdot)^{\frac{1}{p}} (p \in \N_0) $,\/ $\lvert\cdot\rvert$,  $ +, -, \times, /,
\sup(\cdot,\cdot), \inf(\cdot,\cdot) $:
\begin{equation}
\label{eq:f_{sa}}
\fsa^*  :=  \inf_{x \in K} \fsa (x) \enspace.
\end{equation}

Furthermore, when $f = \fpop$ is a multivariate polynomial and $K=\Kpop$ is given by finitely many polynomial inequalities,
Problem~\eqref{eq:f_{sa}} matches the Polynomial Optimization problem (POP):
\begin{equation}
\label{eq:f_{pop}}
\fpop^*  :=  \inf_{x \in \Kpop} \fpop (x) \enspace .
\end{equation}

\subsubsection*{Motivations}
Our ultimate motivation is to automatically verify inequalities
occurring in the proof of Kepler conjecture by Thomas
Hales~\cite{DBLP:journals/dcg/HalesHMNOZ10}.  The formal verification
of Kepler's conjecture is an ambitious goal addressed by the Flyspeck
project~\cite{hales:DSP:2006:432}. Flyspeck is a large-scale
effort needing to tackle various mathematical tools. One particular
difficulty is that Hales' proof relies on hundreds of inequalities, and
checking them requires non-trivial computations. Because of the
limited computing power available inside the proof assistants, it is
essential to devise optimized algorithms that:
\begin{inparaenum}[(1)]
\item verify these inequalities
automatically, and 
\item produce a {\em certificate} for each
inequality, whose checking is computationally reasonably simple. 
\end{inparaenum}


There are numerous other applications to the formal
assessment of such real inequalities; we can point to several other
recent efforts to produce positivity certificates for such problems
which can be checked in proof assistants such as
Coq~\cite{Monniaux_Corbineau_ITP2011}~\cite{Besson:2006:FRA:1789277.1789281},
HOL-light~\cite{harrison-sos} or
MetiTarski~\cite{Akbarpour:2010:MAT:1731476.1731498}.

The Flyspeck inequalities typically involve multivariate polynomials
with some additional transcendental functions; the aim is thus to
compute a lower bound for such expressions. These inequalities
are in general tight, and thus challenging for numerical solvers.
Computing lower bounds in constrained polynomial optimization problems
(POP) is already a difficult problem, which has received much
attention.  Semidefinite relaxation based methods have been developed
in~\cite{DBLP:journals/siamjo/Lasserre01}~\cite{parrilo:polynomials};
they can be applied to the more general class of semialgebraic
problems~\cite{putinar1993positive}.  Alternative approaches are based
on Bernstein polynomials~\cite{Zumkeller:2008:Thesis}.  The task is
obviously more difficult in presence of transcendental functions.
Other methods of choice, not restricted to polynomials, include global optimization by interval methods (see e.g.~\cite{DBLP:journals/rc/Hansen06}), branch and bound methods with Taylor models~\cite{DBLP:journals/mp/CartisGT11}~\cite{Berz:2009:RGS:1577190.1577198}.

In what follows, we will consider the following running example taken
from Hales' proof:


\begin{example}[Lemma$_{9922699028}$ Flyspeck]
\label{ex:9922} Let $K$, $\Delta x$, $l$, and $f$ be defined as follows:
\begin{itemize}
\item $K := [4; 6.3504]^3 \times [6.3504; 8] \times [4; 6.3504]^2$				
\item $\Delta x := x_1 x_4 ( - x_1 +  x_2 +  x_3  - x_4 +  x_5 +  x_6) 
+ x_2 x_5 (x_1  -  x_2 +  x_3 +  x_4  - x_5 +  x_6) 
+ x_3 x_6 (x_1 +  x_2  -  x_3 +  x_4 +  x_5  -  x_6)
- x_2 x_3 x_4  -  x_1 x_3 x_5  -  x_1 x_2 x_6  - x_4 x_5 x_6
$
\item $l(x) := -\frac{\pi}{2} +  1.6294  - 0.2213 (\sqrt{x_2}  +  \sqrt{x_3} +  \sqrt{x_5} +  \sqrt{x_6}  -  8.0) + 0.913 (\sqrt{x_4}  -  2.52) + 0.728 (\sqrt{x_1}  -  2.0)$ 
\item $f (x) :=  l(x) + \arctan \frac{\partial_4 \Delta x }{\sqrt{4 x_1 \Delta x}}$
\end{itemize}
Then, $\forall x \in K,  f(x)  \geq 0.$
\end{example}
\medskip
\if{Note that the inequality above would be
much simpler to check if the map $l$ did not depend on
$x$. Indeed, semialgebraic optimization methods would provide lower
and upper bounds for the argument of $\arctan$, then
we could conclude by monotonicity of $\arctan$ using interval
arithmetic. Here, both $l$ and $t : x \mapsto \arctan \,
\dfrac{\partial_4 \Delta x }{\sqrt{4 x_1 \Delta x}}$ depend on $x$. By
using interval arithmetic addition (without any domain subdivision) on
the sum $l + t$, we obtain a lower bound equal to $-0.87$. Hence, it
is necessary to provide accurate semialgebraic approximations for $t$
and solve one or several instances of Problem~\eqref{eq:f_{sa}}.}\fi

\subsubsection*{Contribution} In this paper, we present a
certification framework, combining Lasserre SDP relaxations
of semialgebraic problems with max-plus approximation
by quadratic functions. 

The idea of max-plus approximation comes
from optimal control: it was originally introduced
by Fleming and McEneaney~\cite{a5}, and developed
by several authors~\cite{a6,a7,curseofdim,PhysRevA.82.042319,conf/cdc/GaubertMQ11}, to represent the value function
by a ``max-plus linear combination'', which
is a supremum of certain basis
functions, like quadratic forms. 
When applied to the present context, this idea leads to
approximate from above and from below
every transcendental function appearing
in the description of the problem
by infima and suprema of finitely
many quadratic forms. 
In that way, we are reduced to a converging sequence
of semialgebraic problems. A geometrical way
to interpret the method is to think of it
in terms of ``quadratic cuts'': quadratic inequalities
are successively added to approximate the graph
of a transcendental function.



The proposed method (Figure~\ref{alg:algo_{samp_optim}})
may be summarized as follows.
Let $f$ be a function and $K$ a box issued from a Flyspeck
inequality, so $f$ belongs to the set of transcendental functions obtained by composition of semialgebraic functions
  with $\arctan$, $\arccos$, $\arcsin$, $\exp$, $\log$, $|\cdot|$,
  $(\cdot)^{\frac{1}{p}} (p \in \N_0) $, $ +, -, \times, /, \sup(\cdot,\cdot), \inf(\cdot,\cdot)$. We alternate steps of approximation, in which an additional quadratic
function is added to the representation, and optimization
steps, in which an SDP relaxation from Lasserre hierarchy
is solved. The information on the location of the optimum
inferred from this relaxation
is then used to refine dynamically the quadratic approximation.
\if{

To solve the resulting Problem~\eqref{eq:f},
  Section~\ref{sect:tr} describes how to underestimate $f$ by a
  semialgebraic function $f_{sa}$ on a compact set $\Ksa$ that
  contains $K$. Then, the feasible points generated by these methods
  are used to refine iteratively the semialgebraic approximations of
  transcendental functions until the needed accuracy is reached. Hence
  it leads to the resolution of a hierarchy of
  Problem~\eqref{eq:f_{sa}} instances. We define this hierarchy by
  adding quadratic constraints to a semialgebraic optimization
  problem. In other words, we optimize an objective function on a
  sequence of semialgebraic relaxations obtained by applying
  successive quadratic cuts on the initial set $K$.
  Based on the properties of these
  relaxations, the following inequality holds:}\fi 
In this way, at each step of the algorithm, 
we refine the following inequalities
	\begin{equation}
	\label{eq:3relax}
		f^* \geq f_{sa}^* \geq f_{pop}^* \enspace,
	\end{equation}
where $f^*$ is the optimal value of the original problem
$f^*_{sa}$ the optimal value of its current semialgebraic approximation,
and $f_{pop}^*$ the optimal value of the SDP relaxation
which we solve. The lower estimate $f_{pop}^*$ does
converge to $f^*$. This follows from a theorem
of Lasserre (convergence of moment SDP relaxations)
and from the consistency of max-plus approximation, see Theorem~\ref{th:transc}.

Max-plus approximation has attracted interest because it may attenuate
the ``curse of dimensionality'' for some structured
problems~\cite{mccomplex}.  Indeed, the estimate
of~\cite{conf/cdc/GaubertMQ11} shows that the number of quadratic
terms needed to reach an $\epsilon$-approximation of a function of $d$
variables is of order $\epsilon^{-d/2}$, where $d$ is the dimension.
Hence, max-plus approximations can be applied to fixed, small
dimensional sub-expressions of complex high dimensional expressions, in
a curse of dimensionality free way. In particular, in the Flyspeck
inequalities involve generally 6 variables, but only univariate
transcendental functions, so $d=1$.

An alternative, more standard approach, is to approximate
transcendental functions by polynomials of a sufficiently high
degree, and to apply SDP relaxations to the polynomial problems
obtained in this way. Further experiments presented in~\cite{victorcicm} indicate that this method is not always scalable. Another alternative approach, which
is quite effective on Flyspeck type inequalities,
is to run branch and bound type algorithms with interval arithmetics.
However, in some instances, this leads to certifying an exponential number
of interval arithmetics computations. Thus, it is of interest
to investigate hybrid methods such as the present one, in
order to obtain more concise certificates.

An important issue, for the practical efficiency of the method,
is the simultaneous tuning of the precision of the max-plus approximation
and of the orders of semidefinite relaxation. 
How to perform optimally this tuning is still not well
understood.  However, we present 
experimental results, both
for some elementary examples as well as non-linear 
        inequalities issued from the Flyspeck project,
giving some indication that certain hard
subclasses of problems (sum of arctan of correlated functions in
many variables) can be solved in a scalable way.

To solve the POP instances, several solvers are available as Gloptipoly~\cite{henrion:hal-00172442} or 
Kojima sparse refinement of the hierarchy of SDP relaxations~\cite{Waki06sumsof},
implemented in the SparsePOP solver~\cite{DBLP:journals/toms/WakiKKMS08}.
These solvers are interfaced with several SDP solvers
  (e.g. SeDuMi~\cite{Sturm98usingsedumi}, CSDP~\cite{Borchers97csdp},
  SDPA~\cite{YaFuNaNaFuKoGo:10}).

The paper is organized as follows. In Section~\ref{sect:pre},
we recall the definition and properties of Lasserre relaxations
of polynomial problems, together with reformulations by Lasserre
and Putinar of semialgebraic problems classes. 
The max-plus approximation, and the main algorithm are presented in Section~\ref{sect:tr}. In Section~\ref{sect:bb}, we show how the algorithm can be combined
with standard domain subdivision methods, to reduce the relaxation gap.
Numerical results are presented in Section~\ref{sect:bench}.

\if{
        $f$ being issued from a particular instance of
        Problem~\eqref{eq:ineq}, we prove that $f_{pop}^* \geq 0$ if
        the relaxation is accurate enough.
 Our complete method described in Section~\ref{sect:bb}
        consists in an iterative decomposition of the set $K$ into
        subsets in which the inequalities to be certified are expected
        to be either tight or coarse. Subsets of $K$ are iteratively
        found by a procedure that checks if a given Taylor
        second-order approximation of the initial function $f$ is
        non-negative on this subsets. In Section~\ref{sect:bench},
        experimental results will illustrate numerical and scalability
        issues. Some elementary examples as well as non trivial
        inequalities issued from the Flyspeck non-linear part are
        presented.
In Section~\ref{sect:pre} we introduce notations and
  preliminaries about Semidefinite Programming (SDP) relaxations to
  explain how to deal with constrained polynomials
  optimization. Section~\ref{sect:sa} reduces the semialgebraic
  optimization Problem~\eqref{eq:f_{sa}} to a constrained polynomial
  optimization problem in a lifted space $K_{pop}$.
  }\fi

\section{NOTATION AND PRELIMINARY RESULTS}
\label{sect:pre}
\if{
We consider the vector space $\mathcal{S}_{n}$ of symmetric $n \times
n$ matrices. It is equipped with the usual inner product $\langle X,
\, Y \rangle = \text{Tr} (XY)$ for $X$, $Y \in \mathcal{S}_{n}$. The
Frobenius norm of a matrix $X \in \mathcal{S}_{n}$ is defined by
$\Vert X \Vert_{F} := \sqrt{\text{Tr} (X^2)}$. A matrix $M \in
\mathcal{S}_{n}$ is called positive semidefinite if $x^T M x \geq 0,
\, \forall x \in \R^n$. In this case, we write $M \succcurlyeq 0$, and
define a partial order by writing $X \succcurlyeq Y$ (resp.\ $X \succ Y$) if and only if $X - Y$ is positive semidefinite (resp.\ positive definite).

Let $\mathcal{B}_d$ denote the basis of monomials for the $d$-degree
real-valued polynomials in $n$ variables :
\begin{equation}
	\label{eq:monbasis}
		1, x_1, x_2, \dots , x_1^2, x_1 x_2, \dots, x_n^2, \dots, x_n^d
	\end{equation}
}\fi
        Let $\R_d[X]$ be the vector space of real forms in $n$
        variables of degree $d$ and $\R[X]$ the set of multivariate
        polynomials in $n$ variables.  
\if{        
        If $\fpop : \R^n \to \R$ is
        a $d$-degree multivariate polynomial, we write
	\begin{equation}
	\label{eq:pdecomp}
	\fpop(x) = \sum_{\alpha} p_{\alpha}x^{\alpha},  
	\end{equation}
	$\text{with} \quad x^{\alpha} := x_1^{\alpha_1} \dots x_n^{\alpha_n} \quad \text{and} \quad \alpha \in \mathcal{F}_d := \{ \alpha \in \N^n : \quad  \sum_{i} \alpha_i \leq d \}$
}\fi
We also define the cone $\Sigma_{d} [X]$ of sums of squares of degree at most $2 d$.
\if{
, i.e.
	\begin{equation}
	\Sigma_{2 d} [X] = \Bigl\{\,\sum_i q_i^2, \, \text{ with } q_i \in \R_d[X] \,\Bigr\}.
	\end{equation}
The set $\Sigma_{2 d} [X]$ is a closed, fully dimensional convex cone in $\R_{2 d}[X]$. We denote by $\Sigma[X]$ the cone of sums of squares of polynomials in $n$ variables.
}\fi
	
\subsection{Constrained Polynomial Optimization Problems and SOS}
\label{sect:pop}\label{sect:sdp}

We consider the general constrained polynomial optimization problem (POP):
\begin{equation}
\label{eq:cons_pop}
\fpop^*  :=  \inf_{x \in \Kpop} \fpop (x),
\end{equation}
where $\fpop : \R^n \to \R$ is a $d$-degree multivariate polynomial,
$\Kpop$ is a compact set defined by polynomials inequalities $g_1(x)
\geq 0,\dots,g_m(x) \geq 0$ with $g_i(x) : \R^n \to \R$ being a
real-valued polynomial of degree $w_i, i = 1,\dots,m$.  We call
$\Kpop$ the feasible set of Problem~\eqref{eq:cons_pop}.
Let $g_0 := 1$. We introduce the $k$-truncated quadratic module $M_k(\Kpop) \subset \R_{2 k}[X]$ associated with $g_1, \cdots, g_m$:
\begin{align*}
M_k (\Kpop) = \Bigl\{\,\sum_{j=0}^m \sigma_j(x) g_j (x): \sigma_j \in \Sigma_{k - \lceil w_j / 2\rceil}[X] \,\Bigr\} 
\end{align*} 

Let $k \geq k_0 := \max( \lceil d / 2
\rceil, \max_{0 \leq j \leq m} \{\lceil w_j / 2\rceil\} )$ and consider the following hierarchy of semidefinite relaxations:
			\[
			Q_k:\left\{			
			\begin{array}{l}
			 \sup\limits_{\mu, \sigma_0, \cdots, \sigma_m} \mu \\			 
			\fpop(x) - \mu  \in M_k(\Kpop),\\
			\end{array} \right.
			\]
and denote by $\sup (Q_k)$ its optimal value.
\begin{theorem}[Lasserre~\cite{DBLP:journals/siamjo/Lasserre01}]
\label{prop:sos}
The sequence of optimal values $(\sup (Q_k))_{k \geq k_0}$ is
non-decreasing. If the quadratic module $M_k(\Kpop)$ is archimedean, then this sequence converges to $\fpop^*$. 
\end{theorem}	
The non-linear inequalities to be proved in the Flyspeck project typically
involve a variable $x$ lying in a box $K \subset \R^n$, thus the archimedean condition holds in our case.
		
\if{
\begin{example}[from $\text{Lemma}_{4717061266}$ Flyspeck]
\label{ex:4717}
Consider the inequality $\forall x \in K, \ \Delta x \geq 0$, where $K := [4,
6.3504]^6$ and $\Delta x$ is the polynomial defined in
Example~\ref{ex:9922}. Using SparsePOP~\cite{DBLP:journals/toms/WakiKKMS08},
the optimal value of $128$ for
the problem $ \inf_{x \in K} \Delta x$ is obtained at
the $Q_2$ relaxation with $\epsilon_{SDP} = 10^{-8}$.
\end{example}
}\fi

\subsection{Semialgebraic Optimization}
\label{sect:sa}

In this section, we recall how the previous approach can be extended to semialgebraic optimization problems by introducing lifting variables.

Given a semialgebraic function $\fsa$, we consider the problem $\fsa^* = \inf_{x \in \Ksa} \fsa (x)$, where $\Ksa := \{x \in \R^n \, : \, g_1(x) \geq
0, \dots, g_m(x) \geq 0 \}$ is a basic semialgebraic set. We suppose that  $\fsa$ is well-defined and thus has a basic semialgebraic lifting. Then, following the approach described in~\cite{DBLP:journals/siamjo/LasserreP10}, we can add auxiliary lifting variables $z_1,\dots,z_p$,
and construct polynomials $ h_1, \dots , h_s \in \R[x,
z_1,\dots,z_p]$ defining the semialgebraic set $\Kpop := \{ (x, \, z_1,\dots,z_p) \in \R^{n+p} : x \in \Ksa, 
 h_1(x, z) \geq 0,\dots, h_s(x, z) \geq 0 \}$, ensuring that $\fpop^* := \inf_{(x, z) \in \Kpop} z_p$ is a lower bound of $\fsa^*$.

\if{			
\begin{definition}[Basic Semialgebraic Lifting]
A semialgebraic function $\fsa$ is said to have a basic semialgebraic lifting if there exist $p, s \in \N$, polynomials $ h_1, \dots , h_s \in \R[X, Z_1,\dots,Z_p]$ and a basic semialgebraic set $\Kpop$ defined by:
\begin{multline*}
\Kpop :=  \{ (x, \, z_1,\dots,z_p) \in \R^{n+p} : x \in \Ksa, \\
 h_1(x, z) \geq 0,\dots, h_s(x, z) \geq 0 \}
\end{multline*}
such that the graph of $\fsa$ (denoted $\Psi_{\fsa}$) satisfies:
\begin{align*}
\Psi_{\fsa} :={} & \{ (x, \, \fsa(x)) : x \in \Ksa\}\\
	  ={} & \{ (x, \, z_p) : (x, \, z) \in \Kpop\} \enspace .
\end{align*} 
\end{definition}

\begin{lemma}[Lasserre, Putinar~\cite{DBLP:journals/siamjo/LasserreP10}]
\label{lemma:bsal}
Let $\mathcal{A}$ be the semialgebraic functions algebra obtained by composition of polynomials with $\lvert\cdot\rvert, +, -, \times, /, \sup(\cdot,\cdot), \inf(\cdot,\cdot), (\cdot)^{\frac{1}{p}} (p \in \N_0)$. Then every well-defined $\fsa \in \mathcal{A}$ has a basic semialgebraic lifting.
\end{lemma}
}\fi

To ensure that the archimedean condition is preserved, we add bound constraints over the lifting variables. These bounds are computed by solving semialgebraic optimization sub-problems. 

\begin{example} [from Lemma$_{9922699028}$ Flyspeck]
\label{ex:atn1}
Continuing Example~\ref{ex:9922}, we consider the function
$\fsa := \frac{\partial_4 \Delta x }{\sqrt{4 x_1 \Delta x}}$ and the set $\Ksa := [4, 6.3504]^3 \times [6.3504, 8] \times [4, 6.3504]^2$. The latter can be equivalently rewritten as 
\[
\Ksa := \{ x \in \R^6 : g_1 \geq 0, \dots, g_{12} \geq 0\}
\]
where $g_1 := x_1 - 4, g_2 := 6.3504 - x_1, \dots , g_{11} := x_6 - 4, g_{12} := 6.3504 - x_6$.

We introduce two lifting variables $z_1$ and $z_2$, respectively representing the terms $\sqrt{4 x_1 \Delta x}$ and $\frac{\partial_4 \Delta x }{\sqrt{4 x_1 \Delta x}}$. We also use a lower bound $m_1$ of $ \inf_{x \in \Ksa} \sqrt{4 x_1 \Delta x}$ and an upper bound $M_1$ of $ \sup_{x \in \Ksa} \sqrt{4 x_1 \Delta x}$ which can be both computed by solving auxiliary sub-problems. 

Now the basic semialgebraic set $\Kpop$ can be defined as follows:

\begin{align*}
\Kpop := {} & \{ (x, z_1 ,z_2) \in \R^{6+2} : 
\begin{multlined}[t]
x \in \Ksa, \, h_l(x, z_1, z_2) \geq 0,\\ 
l = 1,\dots,7 \} 
\end{multlined}
\end{align*}
where the multivariate polynomials $h_l$ are defined by:
\begin{align*}
h_1 & := z_1 - m_1 & h_5 & := z_1\\
h_2 & := M_1 - z_1 & h_6 & := z_2 z_1 - \partial_4 \Delta x \\
h_3 & := z_1^2 - 4 x_1 \Delta x & h_7 & := - z_2 z_1 + \partial_4 \Delta x \\
h_4 & := - z_1^2 + 4 x_1 \Delta x 
\end{align*}
Let $h_0 := 1, \omega_l := \deg h_l, \ 0 \leq l \leq 7$ and define the quadratic module $M_k (\Kpop)$ by: 
\begin{multline*}
M_k (\Kpop) = \Bigl\{\,\sum_{j=1}^{12} \sigma_j(x) g_j (x) + \sum_{l=0}^7 \theta_l(x) h_l (x):  \\
\sigma_j \in \Sigma_{k - 1}[X], 1 \leq j \leq 12, \theta_l \in \Sigma_{k - \lceil \omega_l / 2\rceil}[X], 0 \leq l \leq 7 \,\Bigr\} 
\end{multline*}

Consider the following semidefinite relaxations:
\[
Q_k^{sa}: \left \{
\begin{array}{l}
  \sup\limits_{\mu, \sigma_1, \cdots, \sigma_{12}, \theta_0, \cdots, \theta_7} \mu \\
  z_2 - \mu  \in  M_k (\Kpop)\\
\end{array} \right.
\]
If $k \geq k_0 := \max_{0 \leq l \leq 7} \{\lceil\omega_l / 2\rceil \}= 2$, then as a special case of Theorem~\ref{prop:sos}, the sequence $(\sup (Q_k^{sa}))_{k \geq 2}$ is monotonically non-decreasing and converges to $\fsa^*$. 
A tight lower bound $m_3 = -0.445$ is obtained at the third relaxation.
\end{example}

\section{TRANSCENDENTAL FUNCTIONS UNDERESTIMATORS}
\label{sect:tr}
In this section, we introduce an algorithm allowing to determine that a multivariate transcendental function is positive (Problem~\eqref{eq:ineq}). The algorithm relies on an adaptive basic-semialgebraic relaxation, in which approximations of transcendental
functions by suprema or infima of quadratic forms are iteratively refined.
\if{
using a max-plus semiconvex
approximation. 
In other words, for a given transcendental function
$f$, we build a hierarchy of semialgebraic underestimators with a sup
of quadratic forms.}\fi

\if{
\begin{figure}[t]	
\begin{center}
\begin{tikzpicture}[scale=1]
\draw[->] (-3.5,0) -- (3.5,0) node[right] {$a$};
\draw[->] (0,-1.5) -- (0,2.2) node[above] {$y$};

\draw[color=darkgreen, dotted, thick] plot [domain=-2.5:1.2]  (\x, {0.32*(\x +1)^2 + 1/2* (\x +1) - 0.79 }) node[anchor = north west] {$\parab_{a_1}^+$};
\draw[color=darkgreen, dotted, thick] plot [domain=0.5:3]  (\x, {0.32*(\x -2)^2 + 1/5* (\x -2)  +1.1 }) node[anchor = south west] {$\parab_{a_2}^+$};

\draw[color=whitegreen, dotted, thick] plot [domain=-0.3:3]  (\x, {-0.32*(\x -2)^2 + 1/5* (\x -2) + 1.1 }) node[anchor = north east] {$\parab_{a_2}^-$};
\draw[color=whitegreen, dotted, thick] plot [domain=-1.8:0.5]  (\x, {-0.32*(\x +1)^2 + 1/2* (\x +1)  -0.79 }) node[anchor = north west] {$\parab_{a_1}^-$};

\draw (2,-0.3) node {$a_2$};  \draw (2,3pt) -- (2,-3pt); \draw [black, dashed] (2,0) -- (2,1.1071);
\draw (-1.,-0.3) node {$a_1$};  \draw (-1,3pt) -- (-1,-3pt); \draw [black, dashed] (-1,0) -- (-1,-0.7853);

\draw[color=red!80] plot [domain=-3:3] (\x,{rad(atan(\x))}) node[right] {$\arctan$};
\draw (-3,-0.3) node {$m$}; \draw (3.,-0.3) node {$M$}; \draw (3,3pt) -- (3,-3pt); \draw (-3,3pt) -- (-3,-3pt);
\end{tikzpicture}
\caption{Semialgebraic underestimators and overestimators for $\arctan$}	\label{fig:atn}		
\end{center}
\end{figure}	
}\fi	
	
\subsection{Max-plus Approximation of Semiconvex Functions}
\label{sect:tr1}
Let $\cB$ be a set of functions $\R^n\to \R$, whose elements
will be called \textit{max-plus basis functions}. Given
a function $f:\R^n\to \R$, we look
for a representation of $f$ as a linear combination
of basis functions in the max-plus sense, i.e.,
\begin{equation}
f = \sup_{w\in \cB}(a(w) + w) \label{e-Fenchel}
\end{equation}
where $(a(w))_{w\in \cB}$ is a family of elements of $\R \cup\{-\infty\}$
(the ``coefficients''). 
The correspondence between the function $x \mapsto f(x)$ and the 
coefficient function $w\mapsto a(w)$ is a well studied
problem, which has appeared in various
guises (Moreau conjugacies, generalized Fenchel transforms, Galois correspondences, see~\cite{agk04} for more background).

The idea of max-plus approximation~\cite{a5,mceneaney-livre,a6} is to
choose a space of functions $f$ and a corresponding set $\cB$ of
basis functions $w$, and to approximate from below a given $f$ in this space by a finite max-plus linear combination,
\(
f \simeq  \sup_{w\in \cF}(a(w) + w) 
\)
where $\cF\subset\cB$ is a finite subset. Note that 
$ \sup_{w\in \cF}(a(w) + w) $ is not only an approximation but a valid
lower bound of $f$. 

Following~\cite{a5,a6}, for each constant
$\gamma\in \R$, we shall consider 
the family of quadratic functions $\cB=\{w_{y}\mid y\in \R^n\}$ where 
\[
w_{y}(x):= -\frac{\gamma}{2}\|x-y\|^2  \enspace .
\]
Recall that a function is $\gamma$-semiconvex if and only if the function $x\mapsto \phi(x)+\frac{\gamma}{2}|x|^2$ is convex. Then, it follows from
Legendre-Fenchel duality that the space of functions $f$ which
can be written as~\eqref{e-Fenchel} is precisely the set of lower
semicontinuous $\gamma$-semiconvex functions. 

The transcendental functions
which we consider here are twice
continuously differentiable. Hence, their
restriction to any bounded convex set 
is $\gamma$-semiconvex for a sufficiently large $\gamma$,
so that they can be approximated by finite suprema
of the form $\sup_{w\in \cF}(a(w) + w) $ with $\cF\subset \cB$.
A result of~\cite{conf/cdc/GaubertMQ11} shows that if $N=|\cF|$
basis functions are used, then the best approximation
error is $O(1/N^{2/n})$ (the error is the sup-norm, over any compact set),
provided that the function to be approximated is of class $\mathcal{C}^2$.
Equivalently, the approximation error is of order $O(h^2)$ where $h$
is a space discretization step. Note that the error of max-plus approximation
is of the same order as the one obtained by conventional $P_1$ finite 
elements under the same regularity assumption.
For the applications considered in this paper, $n=1$. 

In this way, starting from a transcendental univariate
elementary function $f \in \mathcal{T}$,  such as $\arctan$, $\exp$, \etc{},
defined on a real bounded interval $I$, we arrive at a 
semialgebraic lower bound of $f$, which is nothing but
a supremum of a finite number of quadratic functions.



\begin{example}
Consider the function $f= \arctan$ on an interval $I := [m, M]$.
For every point $a\in I$, we can find a constant $\gamma$
such that
\[
\arctan (x)  \geq \parab_{a}^-(x):=  -\frac{\gamma}{2} (x-a)^2 +f'(a) (x - a) + f (a)
\enspace .
\]
Choosing $\gamma=\sup_{x\in I} -f''(x)$ always work. However, 
it will be convenient to allow $\gamma$ to depend on the choice of $a$ to get tighter lower bounds. 
Choosing a finite subset $A\subset I$, we arrive at an approximation
\begin{equation}
\label{eq:max_par}
\forall x \in I , \, \arctan \, (a) \geq \max_{a\in A} \, \parab_{a}^-(x) \enspace .
\end{equation} 
\if{
where $\mathcal{C}$ is a set indexing a collection of parabola tangent to the function curve.
Given $i \in \mathcal{C}$, the parabola $\parab_{a_i}^-$ is tangent to the curve at the point $a_i$ and underestimates $f$ on $I$.
More precisely, $\parab_{a_i}^-(x) := c_i (x-a_i)^2 +f'(a_i) (x - a_i) + f (a_i)$, where $f'(a_i) = \frac{1}{1 + a_i^2}$. The coefficient $c_i$ depends on $a_i$ and the curvature variations of $\arctan$ on $I$.

($c$ could be chosen independently of $[m,M]$i.e.\ there exists $c < 0$ such that $x \mapsto \arctan(x) - c x^2$ is convex on $I$. 
}\fi
Semialgebraic overestimators $x \mapsto \min_{a \in A} \parab_{a}^+(x)$ can be defined in a similar way. 
\end{example}

\subsection{An Adaptive Semialgebraic Approximation Algorithm}
\label{sect:alg1}	

We now consider an instance of Problem~\eqref{eq:ineq}. As in Flyspeck inequalities, we assume that $K$ is a box.
We assimilate the objective function $f$ with its abstract syntax tree $t$. We assume that the leaves of $t$ are semialgebraic functions in the set $\mathcal{A}$, and other nodes are univariate transcendental functions ($\arctan$, \etc{}) or basic operations ($+$, $\times$, $-$, $/$). For the sake of the simplicity, we suppose that each univariate transcendental function is monotonic. 

We first introduce the auxiliary algorithm $\tt{samp\_approx}$, presented in Fig.~\ref{alg:algot}. Given an abstract syntax tree $t$ and a box $K$, this algorithm computes lower and upper bounds of $t$ over $K$, and max-plus approximations of $t$ by means of semialgebraic functions. It is also parametrized by a finite sequence of control points used to approximate transcendental functions by means of parabola.

The algorithm $\tt{samp\_approx}$ is defined by induction on the abstract syntax tree $t$. 
When $t$ is reduced to a leaf, \ie{} it represents a semialgebraic function of $\mathcal{A}$, we call the functions $\mathtt{min\_sa}$ and $\mathtt{max\_sa}$ which determine lower and upper bounds using techniques presented in Section~\ref{sect:sa}. In this case, the tree $t$ provides an exact semialgebraic estimator.
If the root of $t$ corresponds to a transcendental function node $r \in \mathcal{T}$ taking a single child $c$ as argument, lower and upper bounds $c_m$
and $c_M$ are recursively obtained, as
well as estimators $c^-$ and $c^+$. Then we apply the function
$\mathtt{build\_par}$ that builds the parabola at the given control points, by using the convexity/semiconvexity properties of $r$,
as explained in Section~\ref{sect:tr1}. An underestimator $t^-$ as well as an
overestimator $t^+$ are determined by composition (so-called
$\mathtt{compose}$ function) of the parabola with $c^-$ and $c^+$. Notice
that the behaviour of $\mathtt{compose}$ depends on the monotonicity properties
of $r$. These approximations $t^-$ and $t^+$ are semialgebraic
functions of $\mathcal{A}$, whence we can also compute their lower and
upper bounds using $\mathtt{min\_sa}$ and $\mathtt{max\_sa}$. 
The last case occurs when the root of $t$ is a binary operation whose arguments are two children $c_1$ and $c_2$. We can apply recursively $\tt{samp\_approx}$ to each child and get semialgebraic underestimators $c_1^-$, $c_2^-$ and
overestimators $c_1^+$, $c_2^+$. Note that when the binary operation is the multiplication or the division, we assume that the estimators of $c_1$ or $c_2$ have a constant sign. We have observed that in practice, all the inequalities that we consider in the Flyspeck project satisfy this restriction.

\begin{figure}[t]
\begin{algorithmic}[1]                    
\Require tree $t$, box $K$, SDP relaxation order $k$, control points sequence $s = x_1,\dots,x_r \in K$
\Ensure  lower bound $m$,  upper bound $M$, lower tree $t^-$, upper tree $t^+$ 
	\If {$t \in \mathcal{A}$}
	\State \Return $\mathtt{min\_sa}$ ($t, \, k$), $\mathtt{max\_sa}$ ($t, \, k$), $t$, $t$
	\ElsIf {$r := \mathtt{root} (t) \in \mathcal{T}$ parent of the single child $c$}					
		\State $m_c$, $M_c$, $c^-$, $c^+ := \mathtt{samp\_approx}(c, K, k, s)$ \label{line:samp_approx}
		\State $\parab^-, \parab^+:= \mathtt{build\_par} (r, m_c, M_c, s)$
		\State $t^-, t^+ := \mathtt{compose}(\parab^-, \parab^+, c^-, c^+)$					
		\State \Return $\mathtt{min\_sa}$ ($t^-, \, k$), $\mathtt{max\_sa}$ ($t^+, \, k$), $t^-$, $t^+$
	\ElsIf {$\mathtt{bop}$ $ := \mathtt{root}$ $(t)$ is a binary operation parent of two children $c_1$ and $c_2$}
		\State $m_{c_i}, M_{c_i}, c_{i}^-, c_{i}^+ := \mathtt{samp\_approx}(c_i, K, k, s)$ for $i \in \{1,2\}$
		\State $t^-, t^+ := \mathtt{compose\_bop} (c_{1}^-, c_{1}^+, c_{2}^-, c_{2}^+)$
		\State \Return $\mathtt{min\_sa}(t^-, k)$, $\mathtt{max\_sa}(t^+, k)$, $t^-$, $t^+$								
	\EndIf
\end{algorithmic}
\caption{$\tt{samp\_approx}$ : recursive semialgebraic max-plus approximation algorithm}   \label{alg:algot}
\end{figure}

Our main optimization algorithm $\tt{samp\_optim}$, presented in Fig.~\ref{alg:algo_{samp_optim}}, relies on $\tt{samp\_approx}$ and chooses the sequence of control points $s$ dynamically. 
At the beginning, the set of control points consists of a single point of the box $K$, chosen so as to minimize the value of the function associated to the tree $t$ among a set of random points (Line~\lineref{line:randeval}). Then, at each iteration of the loop from Lines~\lineref{line:begin_loop} to~\lineref{line:end_loop}, the algorithm $\tt{samp\_approx}$ is called to compute a lower bound $m$ of the function $t$ (Line~\ref{line:samp_approx}). At Line~\lineref{line:argmin}, a minimizer candidate $\Xopt$ of the underestimator tree $t^-$ is computed. It is obtained by projecting a solution $\Xsdp$ of the SDP relaxation of Section~\ref{sect:sa} on the coordinates representing the first order moments, following~\cite[Theorem~4.2]{DBLP:journals/siamjo/Lasserre01}. However, the projection may not belong to $K$ when the relaxation order $k$ is not large enough. This is why tools like SparsePOP use local optimization solver in a post-processing step, providing a point in $K$ which may not be a global minimizer. In any case, $\Xopt$ is then added to the set of control points (Line~\ref{line:add_control_point}). Alternatively, if we are only interested in determining whether the infimum of $t$ over $K$ is non-negative (Problem~\eqref{eq:ineq}), the loop can be stopped as soon as $m \geq 0$.

\begin{figure}[t]                    
\begin{algorithmic}[1]                   
\Require tree $t$, box $K$, $iter_{\max}$ (optional argument)
\Ensure  lower bound $m$, feasible solution $\Xopt$
 	\State $s := [$ $\tt{argmin}$ ($\tt{randeval}$ $t$) $]$ \Comment{$s \in K$} \label{line:randeval}
 	\State $n := 0$
 	\State $m := -\infty$
	\While {$n \leq iter_{\max}$}\label{line:begin_loop}
	\State Choose an SDP relaxation order $k$
	\State $m, \, M, \, t^{-}, \, t^{+} := \tt{samp\_approx}$ $(t, \, K, \, k, \, s)$ 
	\State $\Xopt  := \tt{guess\_argmin}$ $(t^-)$ \Comment{$t^- \, (\Xopt) \simeq m$}\label{line:argmin} 	 		
	\State $s := s \cup \{ \Xopt \}$ \label{line:add_control_point}
	\State $n := n + 1$
	\EndWhile\label{line:end_loop}		
	\State \Return $m, \, \Xopt$
\end{algorithmic}			
\caption{$\tt{samp\_optim}$ : Semialgebraic max-plus optimization algorithm }   \label{alg:algo_{samp_optim}}
\end{figure}

When we call several times $\tt{samp\_approx}$ inside the loop from Lines~\lineref{line:begin_loop} to~\lineref{line:end_loop}, we do not need to always compute recursively the
underestimators and overestimators as well as bounds of all the nodes
and the leaves of the abstract syntax tree. Instead, we ``decorate''
the tree with interval and semialgebraic values containing these
information, based on previous iterations.

\begin{example}[Lemma$_{9922699028}$ Flyspeck]
\label{ex:atn2}
We continue Example~\ref{ex:atn1}. Since we computed lower and upper
bounds ($m$ and $M$) for $\fsa := \frac{\partial_4 \Delta x}{\sqrt{ 4 x_1 \Delta x }}$, we know that the $\fsa$ argument of $\arctan$
lies in $I := [m, M]$. We describe three iterations of the
algorithm. Fig.~\ref{fig:atn3} illustrates the related semialgebraic
underestimators hierarchy.
  
\begin{enumerate}
\setcounter{enumi}{-1} 
\item Multiple evaluations of $f$ return a set of values and we obtain
  a first minimizer guess $x_1 := \tt{argmin}$ ($\tt{randeval}$ $f$) corresponding to the minimal value of the set. \xopt{1}{4.8684}{4.0987}{4.0987} {7.8859}{4.0987}{4.0987}
\item We compute $a_1 := \fsa (x_1) = 0.3962$, get the equation of
  $\parab_1^-$ with $\tt{build_{par}}$ and finally compute
  $ m_1 \leq \min_{x \in K} \{l(x) + \parab_{a_1}^- (\fsa (x)) \}$.
  For $k=2$, we obtain $m_1 =  -0.2816 < 0$ and a new minimizer \xopt{2}{4}{ 6.3504}{6.3504}{6.3504}{6.3504}{6.3504}. 
\item $a_2 := \fsa (x_2) = -0.4449$,  $\parab_{a_2}^-$ and
  $ m_2 \leq \min_{x \in K} \{l(x) + \max_{i \in \{1,
    2\}} \{ \parab$ $_{a_i}^- (\fsa (x))\} \}$. 
  For $k=2$, we get $m_2 = -0.0442 < 0$ and a new minimizer \xopt{3}{4.0121}{ 4.0650}{4.0650}{6.7455}{4.0650}{4.0650}. 
\item $a_3 := \fsa (x_3) = 0.1020$,  $\parab_{a_3}^-$, and
  $ m_3 \leq \min_{x \in K} \{l(x) + \max_{i \in \{1, 2,
    3\}} \{ \parab$ $_{a_i}^- (\fsa (x))\} \}$. 
For $k=2$, we obtain $m_3 = -0.0337 < 0$, obtain a new minimizer $x_4$.
\end{enumerate}

\begin{figure}[t]
\begin{center}
\begin{tikzpicture}[scale=1]	
		
\draw[->] (-3.5,0) -- (3.5,0) node[right] {$a$};
\draw[->] (0,-1.5) -- (0,2.2) node[above] {$y$};

\draw[color=whitegreen, dotted, thick] plot [domain=-0.3:3]  (\x, {-0.32*(\x -2)^2 + 1/5* (\x -2) + 1.1 }) node[anchor = north east] {$\parab_{a_1}^-$};
\draw[color=whitegreen, dotted, thick] plot [domain=-2:0.5]  (\x, {-0.32*(\x +1)^2 + 1/2* (\x +1)  -0.79 }) node[anchor = north west] {$\parab_{a_2}^-$};
\draw[color=whitegreen, dotted, thick] plot [domain=-1:3]  (\x, {-0.32*(\x -1)^2 + 1/2* (\x -1)  +0.79 }) node[anchor = south west] {$\parab_{a_3}^-$};

\draw[color=red!80] plot [domain=-3:3] (\x,{rad(atan(\x))}) node[anchor = south east] {$\arctan$};
\draw (-3,-0.3) node {$m$}; \draw (3.,-0.3) node {$M$}; \draw (3,3pt) -- (3,-3pt); \draw (-3,3pt) -- (-3,-3pt);
\draw (2,-0.3) node {$a_1$};  \draw (2,3pt) -- (2,-3pt); \draw [black, dashed] (2,0) -- (2,1.1071);
\draw (-1.,-0.3) node {$a_2$};  \draw (-1,3pt) -- (-1,-3pt); \draw [black, dashed] (-1,0) -- (-1,-0.7853);
\draw (1.,-0.3) node {$a_3$};  \draw (1,3pt) -- (1,-3pt); \draw [black, dashed] (1,0) -- (1,0.7853);
\end{tikzpicture}
\caption{A hierarchy of Semialgebraic Underestimators for $\arctan$}	\label{fig:atn3}
\end{center}
\end{figure}
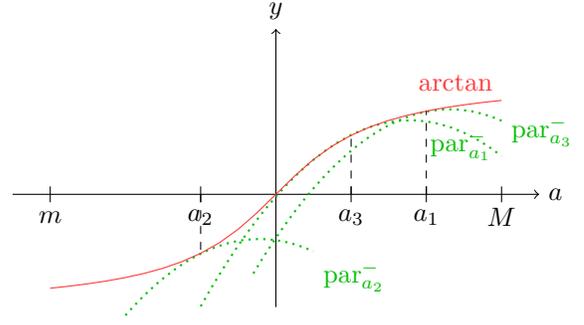
\end{example}		

%


We denote by $t^-_i$ the underestimator computed at the $i^\text{th}$ iteration of the algorithm $\tt{samp\_optim}$, and by $\Xopt^i$ the corresponding minimizer candidate. 
\begin{theorem}[Convergence of $\tt{samp\_optim}$]
\label{th:transc}
Assume that at each iteration $i$, the SDP relaxation order is chosen
to be large enough so that $\Xopt^i$ is a global minimizer of $t^-_i$. 
Then every accumulation point of the sequence of $(\Xopt^i)_i$ is a global minimizer of $t$ over~$K$.
\end{theorem}

Theorem~\ref{th:transc} can be proved using the convergence of Lasserre hierarchy of SDP relaxation~\cite{DBLP:journals/siamjo/Lasserre01}, together with the following lemma:
\begin{lemma}
\label{lemma:space1}
There exists a positive constant $C$ such that,
\begin{equation}\label{eq:space1}
\forall i \geq 1, \forall x \in K, \ t(x) - t_i^-(x) \leq C d(x, S_i)^2,
\end{equation}
where $d$ denotes the euclidean distance and $S_i$ is the set of points $\Xopt^1,\dots,\Xopt^i$.
\end{lemma}
	
The time complexity of our algorithm strongly depends on the relaxation order $k$. Indeed, if $p$ is the number of the control points, then the number of moment variables in the SDP problem $Q_k$ is in $O((n + p)^{2 k})$, and the size of linear matrix inequalities involved are in $O((n + p)^k)$.
The complexity of $\mathtt{samp\_optim}$ is therefore exponential in~$k$. 
%
Notice that there are several ways to decrease the size of these matrices. First, symmetries in SDP relaxations for polynomial
optimization problems can be exploited to replace one SDP problem $Q_k$ by
several smaller SDPs~\cite{DBLP:journals/corr/abs-1103-0486}. Notice it is possible only if the multivariate polynomials of the initial problem are invariant under the action of a finite subgroup $G$ of the group $GL_{n+p}(\R)$. Furthermore, one can exploit the structured sparsity of the
problem to replace one SDP
problem $Q_k$ by an SDP problem of
size $O (\kappa^ {2 k})$ where $\kappa$ is the average size
of the maximal cliques correlation pattern of the polynomial
variables (see~\cite{DBLP:journals/toms/WakiKKMS08}). 
\section{REFINING BOUNDS BY DOMAIN SUBDIVISION}
\label{sect:bb}


A small relaxation order ensures fast computation of the lower bounds but the relaxation gap may remain too high to ensure
the convergence of the algorithm. This is particularly critical when we want to certify that a given transcendental multivariate function is non-negative. In this section, we explain how to reduce the relaxation gap using domain subdivision in order to solve problems of the form~\eqref{eq:ineq}.


Suppose that the algorithm $\tt{samp\_optim}$ returns a negative lower bound $m$ and a global minimizer candidate $x_c^*$. Our approach consists in cutting the initial box $K$ in several boxes $(K_i)_{1 \leq i \leq c}$. We explain the partitioning of $K$ with the following heuristic. 

Let $\mathcal{B}_{x_c^*,\, r}$ be the intersection of the 
$L_{\infty}$-ball of center $x_c^*$ and radius $r$ with the set $K$. Then, let $f_{x_c^*, r}$ be the quadratic form defined by:
\begin{align*}
    f_{x_c^*,\, r} : \mathcal{B}_{x_c^*,\, r} \longrightarrow{} & {} \R \\
    x \longmapsto{} & {} f(x_c^*) + \Df (x_c^*) \, (x - x_c^*) \\
    & + \frac{1}{2}(x - x_c^*)^{T} \Hf (x_c^*) \, (x - x_c^*) \\
    & + \frac{1}{2} \lambda(x - x_c^*)^2  
\end{align*}
with $\lambda$ given by:
\begin{equation}
\label{eq:lambda}
\lambda := \min_{x \in \mathcal{B}_{x_c^*, \, r}} \{ {\lambda_{\min} ( \Hf (x) - \Hf (x_c^*)  )} \}
\end{equation}

\begin{lemma}
\label{th:tm}
$\forall x \in \mathcal{B}_{x_c^*,\, r}, \ f(x) \geq f_{x_c^*,\, r}$.
\end{lemma}



To underestimate the value of $\lambda$, we determine an interval matrix
$\widetilde{\Hf} := ([\underline{d_{i j}}, 
\overline{d_{i j}} ])_{1 \leq i, j \leq n}$ containing coarse bounds of the Hessian difference
$( \Hf (x) - \Hf (x_c^*) )$ on
$\mathcal{B}_{x_c^*, \, r}$ using interval arithmetic or $\tt{samp\_approx}$ with a small number of control points and a low SDP relaxation order. We then apply on $\widetilde{\Hf}$ a robust SDP method on interval matrix described by
Calafiore and Dabbene in~\cite{springerlink:10.1007/s10957-008-9423-1}, and obtain a lower bound $\lambda'$ of $\lambda$.
%

%

By dichotomy and using Lemma~\ref{th:tm}, we can finally compute the $L_{\infty}$-ball $\mathcal{B}_{x_c^*, \, r}$ of maximal radius $r$ such that the underestimator $f_{x_c^*,\, r}$ is non-negative on $\mathcal{B}_{x_c^*, \, r}$.
\if{
An illustration of our subdivision algorithm is given in
Fig.~\ref{fig:bb2d} in the two dimensional case.

\begin{figure}[t]    
\begin{center}
\begin{tabular}{ccc}
\begin{tikzpicture}[scale=3, baseline={($(current bounding box.center)+(0,-1ex)$)}]							
\draw [black, pattern=checkerboard light gray] ( (0,0) -- (1,0) -- (1,1) -- (0,1) -- cycle;
\draw [black, fill=white] (0.2,0.2) -- (0.7,0.2) -- (0.7,0.7) -- (0.2,0.7) -- cycle;
\draw (0.3,0.3) node[above right] {$x_c^*$} node{$\bullet$};
\draw (0.7,0.7) node[below left] {$\mathcal{B}_{x_c^*,\, r}$};				
\end{tikzpicture}& $\Rightarrow$&
\begin{tikzpicture}[scale=3, baseline={($(current bounding box.center)+(0,-1ex)$)}]							
\draw [black,fill=darkred,  pattern=checkerboard light gray] (0,0) -- (1,0) -- (1,1) -- (0,1) -- cycle;
\draw [black,fill=white] (0.2,0.2) -- (0.7,0.2) -- (0.7,0.7) -- (0.2,0.7) -- cycle;
\draw (0.3,0.3) node[above right] {$x_c^*$} node{$\bullet$};
\draw (0.5,0.7) node[below] {$K_0$};				
\draw [black, dashed] (0.2,1) -- (0.2,0.7);
\draw [black, dashed] (0.2,0) -- (0.2,0.2);
\draw [black, dashed] (0.7,1) -- (0.7,0.7);
\draw [black, dashed] (0.7,0) -- (0.7,0.2);
\draw (0.1,0.5) node {$K_1$};\draw (0.45,0.1) node {$K_2$};\draw (0.45,0.85) node {$K_3$};\draw (0.85,0.5) node {$K_4$};                       
\end{tikzpicture}
\end{tabular}
\end{center}
\caption{A two dimensional example for our box subdivision} \label{fig:bb2d}
\end{figure}
}\fi
\section{RESULTS}
\label{sect:bench}

We next present the numerical results obtained with our method for both small and medium-sized inequalities taken from the Flyspeck project. 

In Tables~\ref{table:atn1} and~\ref{table:atn3}, the inequalities are indexed by the first four digits of the hash code. We also indicate in subscript the number of variables involved in each inequality. The integer $n_{\mathcal{T}}$ represents the number of transcendental univariate nodes in the corresponding abstract syntax trees. The parameter $k_{\max}$ is the highest SDP relaxation order used to 
solve the polynomial optimization problems with SparsePOP. We denote by $n_{\mathit{pop}}$ the total number of POP that have to be solved to prove the
inequality, and by $n_{\mathit{cuts}}$ the number of domain cuts that are performed during the subdivision algorithm. Finally, $m$ is the lower bound of
the function $f$ on $K$ that we obtain with our method, i.e.\ the minimum of all the computed lower bounds of $f$ among the $n_{\mathit{cut}}$ sub-boxes of $K$.

The inequalities reported in Table~\ref{table:atn1} are similar to the one presented in Example~\ref{ex:9922}. They all consist in the addition of the function $x \mapsto \arctan \frac{\partial_4 \Delta x }{\sqrt{4 x_1 \Delta x}}$ with an affine function over $\sqrt{x_i}$ ($1 \leq i \leq 6$). 
%
\begin{table}[!ht]
\begin{center}
\caption{Results for small-sized Flyspeck inequalities}
\begin{tabular}{|l|c|c|c|c|c|c|}
\hline
Ineq. id & $n_{\mathcal{T}}$  &$k_{\max}$  &$n_{\mathit{pop}}$  & $n_{\mathit{cuts}}$ &  $m$ & time\\
\hline
\hline
$9922_6$ & $1$  & $2$ & $222$ & $27$   &  $3.07 \times 10^{-5}$ & $20 \, min$ \\
\hline			
$3526_6$ & $1 $ & $2$ & $156$ & $17$  & $4.89 \times 10^{-6}$ & $13 \, min$ \\
\hline
$6836_6$ & $1$ & $2$ & $173$ & $22$  & $4.68 \times 10^{-5}$ & $14 \, min$ \\
\hline
$6619_6$ & $1$  & $2$ & $163$ & $21$ & $4.57 \times 10^{-5}$ & $13 \, min$ \\
\hline
$3872_6$ & $1$  & $2$ & $250$ & $30$  & $7.72 \times 10^{-5}$ & $20 \, min$ \\
\hline
$3139_6$ & $1$  & $2$ & $162$ & $17$ & $1.03 \times 10^{-5}$ & $13 \, min$ \\
\hline
$4841_6$ & $1$  & $2$ & $624$ & $73$ & $2.34 \times 10^{-6}$ & $50 \, min$ \\
\hline
$3020_5$ & $1$  & $3$ & $80$ & $9$ & $2.96 \times 10^{-5}$ & $31 \, min$ \\
\hline
$3318_6$ & $1$  & $3$ & $26$ & $2$  & $3.12 \times 10^{-5}$ & $1.2 \, h$ \\
\hline
\end{tabular}
\label{table:atn1}
\end{center}			
\end{table}	

Table~\ref{table:atn3} provides the numerical results obtained on medium-sized Flyspeck inequalities. Inequalities $7394_i$ ($3 \leq i \leq 5$) are obtained from a same inequality $7394_6$ involving six variables, by instantiating some of the variables by a constant value. Inequalities $7726_6$ and $7394_6$ are both of the form $l(x) + \sum_{i = 1}^3 \arctan(q_i(x))$ where $l$ is an affine function over $\sqrt{x_i}$, where $q_1(x) := \frac{\partial_4 \Delta x }{\sqrt{4 x_1 \Delta x}}$, $q_2(x) := q_1(x_2, x_1, x_3, x_5, x_4, x_6)$, and $q_3(x) := q_1(x_3, x_1, x_2, x_6, x_4, x_5)$. 

\begin{table}[!ht]
\begin{center}
\caption{Results for medium-size Flyspeck inequalities}
\begin{tabular}{|l|c|c|c|c|c|c|}
\hline
Ineq. id & $n_{\mathcal{T}}$  &$k_{\max}$  &$n_{pop}$  & $n_{cuts}$	 &  $m$ & time\\
\hline
\hline			
 $7726_6$ & $3$  & $2$ & $450$ & $70$   &  $1.22 \times 10^{-6}$ & $3.4 \, h$ \\
\hline
 $7394_3$ & $3$  & $3$ & $1$ & $0$   &  $3.44 \times 10^{-5}$ & $11 \, s$ \\
\hline			
$7394_4$ & $3 $ & $3$ & $47 $ & $10$  & $3.55 \times 10^{-5}$ & $26 \, \textit{min}$ \\
\hline
 $7394_5$ & $3$  & $3$ & $290$ & $55$  & $3.55 \times 10^{-5}$ & $12 \, h$\\
\hline						
\end{tabular}	
\label{table:atn3}		
\end{center}	
\end{table}

\begin{table}[!ht]
\begin{center}
\caption{Comparison results for random examples}
\begin{tabular}{|c||c|c||c|c|}
\hline
\multirow{2}{*}{$n$}
 & 
\multicolumn{2}{c||}{$\tt{samp\_approx}$ with $k = 3$}&
\multicolumn{2}{c|}{$\tt{intsolver}$}
\\
\cline{2-5}
& $m$  &time  &  $m$ & time\\
\hline
\hline            
 $3$ &  $0.4581 $ & $3.8 \, s$ &  $0.4581 $ & $15.5 \, s $ \\
\hline            
 $4$ & $0.4157 $ & $12.9 \, s $   &  $0.4157 $ & $172.1 \, s $ \\
\hline            
 $5$ & $0.4746 $ & $1 \, \textit{min} $    &  $0.4746 $ & $10.2 \, \textit{min} $ \\
\hline            
 $6$ &  $0.4476 $ & $4.6 \, \textit{min} $    &  $ 0.4476$ & $3.4\, h $ \\
\hline
\end{tabular}
\label{table:rnd}
\end{center}
\end{table}
	
In Table~\ref{table:rnd}, we compared our algorithm with the MATLAB toolbox $\tt{intsolver}$~\cite{Montanher09} (based on the Newton interval method~\cite{Hansen198389}) for random inequalities involving two transcendental functions. We denote by $n$ the number of variables, and by $m$ the lower bound that we obtain. The functions that we consider are of the form $x \mapsto \arctan(p(x)) + \arctan(q(x))$, where $p$ is a four-degree polynomial and $q$ is a quadratic form. All variables lie in $[0, 1]$. Both $p$ and $q$ have random coefficients (taken in $[0,1]$) and are sparse. The results indicate that for such examples, our method may outperform interval arithmetic.



\section{CONCLUSION}

We proposed a hybrid method to certify tight non-linear inequalities, 
combining SDP and approximation of semiconvex functions by suprema of quadratic forms (max-plus basis method, originating from optimal control).
The proposed approach
bears some similarity with the ``cutting planes'' proofs in 
combinatorial optimization, the cutting planes being now
replaced by quadratic inequalities. 
This allowed us to solve both small and intermediate size
inequalities of the Flyspeck project, with 
a moderate order of SDP relaxation.

\if{
Our subdivision algorithm performs an iterative refinement of the
initial domain where the inequality is either tight or coarse. Hence,
the parts of the domain where it is difficult to solve are focused by
the iterative partitioning. The initial objective function is
underestimated by quadratic forms, using Taylor Second Order
approximations.}\fi
Several improvements are possible. The approximation technique used 
here could be also applied recursively to some semialgebraic subexpressions,
in order to reduce the degree of the POP instances. 
\if{
Furthermore, a
well-suited partitioning of the domain will be done, using also an
iterative scheme. Some techniques consisting in sphere or ellipsoids
covers of the domain have been already
investigated~\cite{citeulike:10833589}~\cite{DBLP:journals/mp/CartisGT11}.
}\fi
\if{
Other Branch and Bound algorithms could be implemented. For instance,
the well-known Moore algorithm~\cite{Moore1962IAa} mainly used to refine a
domain after studying the sign of the partial derivatives.
}\fi 
\if{
Moreover, randomization techniques as well as the so-called
``joint+marginal'' approach~\cite{DBLP:conf/cdc/LasserreT10} could be
implemented in order to detect better feasible points for the
optimization problems.
}\fi

Finally, we plan to interface the present framework with the Coq proof assistant~\cite{CoqProofAssistant}, by exploiting formally certified symbolic-numeric algorithms~\cite{KLYZ09}.
We believe that
\if{We should finally mention that some recent discussions with members of the
Flyspeck project let us think that }\fi
hybrid certification techniques (combinations of formal methods)
could be a suitable option to formalize the remaining non-linear
inequalities within reasonable amount of CPU time in order to complete
the remaining part of the formal verification of the proof of Kepler conjecture.

\label{sect:concl}

\section*{Acknowledgements}
The authors thank the anonymous referees for helpful comments and suggestions to improve this paper.

\bibliography{mybib}

\end{document}